\documentclass{svproc}
\usepackage{url}

\begin{document}
\mainmatter              % start of a contribution
\title{Towards discrete octonionic analysis}
\titlerunning{Towards discrete octonionic analysis}  % abbreviated title (for running head)
\author{Rolf S\"oren Krau\ss har\inst{1} \and Anastasiia Legatiuk\inst{2} \and
Dmitrii Legatiuk\inst{1}}
\authorrunning{R.S. Krau\ss har, A. Legatiuk, D. Legatiuk} % abbreviated author list (for running head)
\tocauthor{Rolf S\"oren Krau\ss har, Anastasiia Legatiuk, Dmitrii Legatiuk}
\institute{Universit\"at Erfurt, Chair of Mathematics, Erfurt, Germany\\
\and
Bauhaus-Universit\"at Weimar, Chair of Applied Mathematics, Weimar, Germany}

\maketitle              % typeset the title of the contribution

\begin{abstract}
In recent years, there is a growing interest in the studying octonions, which are 8-dimensional hypercomplex numbers forming the biggest normed division algebras over the real numbers. In particular, various tools of the classical complex function theory have been extended to the octonionic setting in recent years. However not so many results related to a discrete octonionic analysis, which is relevant for various applications in quantum mechanics, have been presented so far. Therefore, in this paper, we present first ideas towards discrete octonionic analysis. In particular, we discuss several approaches to a discretisation of octonionic analysis and present several discrete octonionic Stokes' formulae.
% We would like to encourage you to list your keywords within
% the abstract section using the \keywords{...} command.
\keywords{octonions, discrete Clifford analysis, discrete operators, Stokes' formula, discrete octonions}
\end{abstract}

\section{Introduction}
As very well-known, complex analysis provides a very powerful toolkit to study numerous boundary value problems arising in classical harmonic analysis in the two-dimensional case. Motivated by modern problems of engineering and physics, there has been a rapidly growing interest in developing higher-dimensional versions of complex function theory to extend the classical toolkit for a successful treatment of higher-dimensional problems. While engineering mainly focusses on three-dimensional settings, modern physics, for example particle physics, also require tools in the context of dimensions $n > 3$. Einstein's relativity theory already requires four dimensions, including time. Also the standard model of particle physics of electro-weak action requires four dimensions. An enormous challenge in modern physics however, is to understand how gravity can be incorporated on the level of particle physics. Studying of this problems leads to the consideration of even higher dimensional settings, such as string and super-string theory, where the latter requires $12$ dimensions. More recent models of a generalised standard model give stronger indications to use an eight dimensional model, see for example works \cite{Burdik,G,NSM} which shall provide one motivation for this paper from the physical point of view.\par  
From the mathematical point of view, there are several possibilities to extend complex numbers and complex analysis to higher dimensions. One approach is to work with associative Clifford algebras leading to several function theories that consider functions defined on open subsets of an arbitrary dimensional vector space $\bbbr^{n+1}$ that take values in a $2^n$-dimensional Clifford algebra $\bbbr_n$. On this way, a higher-dimensional version of the Cauchy-Riemann operator $D := \sum\limits_{i=0}^n \partial_{x_i} e_i$ factorising the $n+1$-dimensional Laplacian in the form of an elliptic first order differential operator is considered. Its function theory is widely known as {\itshape Clifford analysis}, see for instance \cite{BDS}, and offers many powerful generalisations of complex function theory, such as a Cauchy integral formula, Taylor and Laurent series expansions, a residue theory and a toolkit to study operators of Calderon-Zygmund type on strongly Lipschitz surfaces. A series of textbooks, see for example \cite{GS2}, presents a toolkit of related integral operators that can be used to tackle associated boundary value problems. Recently a lot of progress has been made in also elaborating discrete versions of Clifford analysis which also opened the door to apply these function theoretic tool numerically in bounded and unbounded domains, see \cite{Brackx_1,Cerejeiras_1,CKKS,Cerejeiras_2,Cerejeiras,FGHK06,Faustino_1,FKS07,Guerlebeck_1} among others.\par
However, besides the use of associative Clifford algebras, there are also other possibilities of generalising complex function theory to higher dimensions. If the Cayley-Dickson duplication process to the complex numbers is applied, then we first arrive at the four-dimensional Hamiltonian quaternions, which, however, still is a Clifford algebra; and after applying it once more, we obtain a new algebra, namely the {\itshape octonions}, see \cite{Baez}. Octonions are not any more associative -- so, they are neither a Clifford algebra nor representable with matrices. However, they still form a normed non-associative division algebra having no zero-divisors. From the recent viewpoint of generalised particle physics, see again for example \cite{Burdik,G,NSM}, octonions seem to offer a more adequate model for a unified description of particle physics including gravity, see also \cite{GTBook}. However, there is still a lack of results on the level of octonionic function theory.\par
According to our knowledge, the first contribution to introduce an octonionic generalisation of complex function theory was provided by P. Dentoni and M. Sce in 1973 in \cite{DS}, where a Cauchy integral formula for null-solutions to the octonionic Cauchy-Riemann operator has been presented. Later, a lot of fundamental contributions were provided by K. Nono \cite{Nono} in 1988, and the school of Xingmin-Li, Li-Zhong Peng and their co-authors starting with 2000 up to now, see \cite{XL2000,XL2001,XL2002,XZL}. In these papers, for instance, generalisations of a Cauchy integral formula together with Plemelj projection formulas and with some basic applications to Calderon-Zygmund type operators \cite{XLT2008} including a generalisation of the three-line theorem from J. Peetre \cite{XL2004}, as well as Taylor and Laurent series expansions have intensively been studied \cite{XL2001}. More recently, J. Kauhanen and H.Orelma started to look more intensively at some elementary octonionic boundary value problems and analysed more precisely the algebraic structure of the set of octonionic null-solutions of Cauchy-Riemann operators, see for example \cite{Kauhanen_1,Kauhanen_2,Kauhanen_3}. As also mentioned by J. Kauhanen and H. Orelma, in contrast to Clifford analysis, octonionic monogenic functions do not form $\bbbo$-modules but only $\bbbr$-modules. This fact has a strong influence on the study of generalised Hilbert spaces in the octonionic setting, which is a topic of very recent research, see for example \cite{ConKra2021,FL,QR2022}. For solving related boundary values in practice, it is necessary to apply discretised versions of the related octonionic operators.\par
Although a discretisation of octonionic analysis is important for practical use of function theoretic tools, to the best of our knowledge, no results related to a discrete octonionic analysis have been presented so far. Therefore, the aim of this short paper is presenting some first results in this direction. In particular, we introduce discretised versions of the octonionic Cauchy-Riemann operators and establish a generalised version of the Stokes' formula. As it will be clearly seen, already at this level, we encounter substantial differences to the classical discrete Clifford analysis: because of the non-associativity of octonionic multiplication, discretisation of octonionic analysis needs to be discussed more carefully. Additionally, we will also indicate the difference to the continuous case, which appear due to working with forward and backward Cauchy-Riemann operators. Thus, this paper serves as a first step for developing discrete octonionic analysis, and results presented here will be further extended in future work.\par

\section{Preliminaries and notations}

\subsection{Continuous octonionic analysis}

Before introducing discrete constructions, let us briefly recall some basic information about octonions $\bbbo$ and continuous octonionic analysis. Let us consider $8$-dimensional Euclidean space $\bbbr^{8}$ with the basis unit vectors $e_{k}$, $k=1,2,\ldots,8$ and points $\mathbf{x}=(x_{1}, x_{2},\ldots, x_{8})$. Then in real coordinates, octonions are expressed in the form
\begin{eqnarray*}
x = x_{0}e_{0}+x_{1}e_{1}+x_{2}e_{2}+x_{3}e_{3}+x_{4}e_{4}+x_{5}e_{5}+x_{6}e_{6}+x_{7}e_{7},
\end{eqnarray*}
where $e_{4}=e_{1}e_{2}$, $e_{5}=e_{1}e_{3}$, $e_{6}=e_{2}e_{3}$ and $e_{7}=e_{4}e_{3}=(e_{1}e_{2})e_{3}$. Additionally we have $e_{i}^{2}=-1$ and  $e_{0}e_{i}=e_{i}e_{0}$ for all $i=1,\ldots,7$, and $e_{i}e_{j}=-e_{j}e_{i}$ for all mutual distinct $i,j\in\left\{1,\ldots,7\right\}$. Table~\ref{Table_octonions} shows multiplication rules for real octonions. As it can be clearly seen from this table, multiplication of octonions is not associative, precisely we have $(e_{i}e_{j})e_{k}=-e_{i}(e_{j}e_{k})$.\par
\begin{table}
\caption{Multiplication table for real octonions $\bbbo$}
\label{Table_octonions}
\begin{center}
\begin{tabular}{|c|cccccccc|}
\hline
$\cdot$ & $e_{0}$ & $e_{1}$ & $e_{2}$ & $e_{3}$ & $e_{4}$ & $e_{5}$ & $e_{6}$ & $e_{7}$ \\[2pt]
\hline
$e_{0}$ & $1$ & $e_{1}$ & $e_{2}$ & $e_{3}$ & $e_{4}$ & $e_{5}$ & $e_{6}$ & $e_{7}$ \\
$e_{1}$ & $e_{1}$ & $-1$ & $e_{4}$ & $e_{5}$ & $-e_{2}$ & $-e_{3}$ & $-e_{7}$ & $e_{6}$ \\
$e_{2}$ & $e_{2}$ & $-e_{4}$ & $-1$ & $e_{6}$ & $e_{1}$ & $e_{7}$ & $-e_{3}$ & $-e_{5}$ \\
$e_{3}$ & $e_{3}$ & $-e_{5}$ & $-e_{6}$ & $-1$ & $-e_{7}$ & $e_{1}$ & $e_{2}$ & $e_{4}$ \\
$e_{4}$ & $e_{4}$ & $e_{2}$ & $-e_{1}$ & $e_{7}$ & $-1$ & $-e_{6}$ & $e_{5}$ & $-e_{3}$ \\
$e_{5}$ & $e_{5}$ & $e_{3}$ & $-e_{7}$ & $-e_{1}$ & $e_{6}$ & $-1$ & $-e_{4}$ & $e_{2}$ \\
$e_{6}$ & $e_{6}$ & $e_{7}$ & $e_{3}$ & $-e_{2}$ & $-e_{5}$ & $e_{4}$ & $-1$ & $-e_{1}$ \\
$e_{7}$ & $e_{7}$ & $-e_{6}$ & $e_{5}$ & $-e_{4}$ & $e_{3}$ & $-e_{2}$ & $e_{1}$ & $-1$ \\[2pt]
\hline
\end{tabular}
\end{center}
\end{table}\par
There are several possibilities to extend the classical function theory to octonions. One way consists of the Riemann-approach, following the line of investigation of P. Dentoni and M. Sce \cite{DS}, K. Nono \cite{Nono}, the school of Xingmin-Li and Zhong Peng, see for instance \cite{XL2000,XL2001} and others. In their spirit one may introduce
\begin{definition}[Octonionic monogenicity]
Let $U \subseteq \bbbo$ be open. A function $f:U \to \bbbo$ is called left (right) octonionic monogenic if ${\cal{D}} f = 0$ (esp. $f {\cal{D}} = 0$). Here, 	${\cal{D}}:= \frac{\partial }{\partial x_0} + \sum\limits_{i=1}^7 e_i \frac{\partial }{\partial x_i}$ is the octonionic first order Cauchy-Riemann operator. If $f$ satisfies $\overline{{\cal{D}}}f = 0$ (resp. $f\overline{\cal{D}} = 0$), then we call $f$ left (right) octonionic anti-monogenic.
\end{definition}\par
In contrast to Clifford analysis, where one considers null-solutions to the Cauchy-Riemann operator defined on the paravector space $\bbbr\oplus\bbbr^7$ with values in the Clifford algebra $\mathcal{C}\ell_7$, which is a real vector spaces isomorphic to $\bbbr^{128}$, the octonionic approach really considers maps from $\bbbr^8$ into $\bbbr^8$. Another essential difference is the fact that left (right) octonionic monogenic functions do neither form a right nor a left $\bbbo$-module. Following for instance J. Kauhanen and H. Orelma in \cite{Kauhanen_3}, one can take as a very simple counterexample: the function $f(x):= x_1 - x_2 e_4$. Then we have ${\cal{D}}[f(x)] = e_1 - e_2 e_4 = e_1 - e_1 = 0$. However, $g(x):=(f(x))\cdot e_3 = (x_1 - x_2 e_4) e_3 = x_1 e_3 - x_2 e_7$ satisfies ${\cal{D}}[g(x)] = e_1 e_3 - e_2 e_7 = e_5 -(-e_5) = 2 e_5 \neq 0$.\par
As already mentioned in the classical paper \cite{XL2000}, the lack of associativity prevents us from getting a direct analogue of Stokes' formula in the octonionic setting. Even if both ${\cal{D}} f = 0$ and $g {\cal{D}} = 0$, we do not have in general 
\begin{eqnarray*}
\int\limits_{\partial G} g(x) \; (d\sigma(x) f(x)) = 0 \mbox{ nor } \int\limits_{\partial G} (g(x) d\sigma(x)) \; f(x) = 0.
\end{eqnarray*}
Instead, quoting from \cite{XLT2008}, we obtain the following relation
\begin{equation}
\label{Stokes_formula_continuous}
\int\limits_{\partial G} g(x) \; (d\sigma(x)  f(x)) = \int\limits_G \Bigg(   
g(x)({\cal{D}} f(x)) + (g(x){\cal{D}})f(x)  - \sum\limits_{j=0}^7 [e_j, {\cal{D}}g_j(x),f(x)]
\Bigg) dV,
\end{equation}
where $[a,b,c] := (ab)c - a(bc)$ is the so-called {\itshape associator} (which would vanish in the cases of associativity). Although the associator appears in most of octonionic constructions, it is nonetheless possible to introduce specific structures, where the associator would vanish. For example, it has been pointed out in \cite{XL2000}, that considering the two functions being octonionic monogenic and Stein-Weiss conjugate harmonics, i.e. $\frac{\partial g_j}{\partial x_i} = \frac{\partial g_i}{\partial x_j}$ for all $0\leq i<j\leq 7$, the associator will vanish.\par
Moreover, it is still possible to obtain a generalisation of the Cauchy's integral formula to octonionic setting \cite{Nono,XL2002}:  
\begin{proposition}[Cauchy's integral formula]\label{cauchy1}
Let $U \subseteq \bbbo$ be open and $G \subseteq U$ be an $8$-D compact oriented manifold with a strongly Lipschitz boundary $\partial G$. If $f: U \to \bbbo$ is left octonionic monogenic, then for all $x \in G$
\begin{eqnarray*}
f(x)= \frac{3}{\pi^4} \int\limits_{\partial G} q_{\bf 0}(y-x) \Big(d\sigma(y) f(y)\Big).
\end{eqnarray*}
\end{proposition}
However, we have to emphasise carefully on the fact that putting the parenthesis differently, leads to the different formula  
\begin{eqnarray*}
\frac{3}{\pi^4} \int\limits_{\partial G} \Big( q_{\bf 0}(y-x) d\sigma(y)\Big) f(y)  =  f(x) +  \int\limits_G \sum\limits_{i=0}^7 \Big[q_{\bf 0}(y-x),{\cal{D}}f_i(y),e_i  \Big] dy_0 \cdots dy_7, 
\end{eqnarray*}
involving the associator again.\par

\subsection{Discretisation of octonionic analysis}

Let us now introduce a discrete setting for octonions. Consider the unbounded uniform lattice $h\bbbz^{8}$ with the lattice constant $h>0$, which is defined in the classical way as follows
\begin{eqnarray*}
h \bbbz^{8} :=\left\{\mathbf{x} \in {\bbbr}^{8}\,|\, \mathbf{x} = (m_{1}h, m_{2}h,\ldots, m_{8}h), m_{j} \in \bbbz, j=1,2,\ldots,8\right\}.
\end{eqnarray*}
Next, we define the classical forward and backward differences $\partial_{h}^{\pm j}$ as
\begin{equation}
\label{Finite_differences}
\begin{array}{lcl}
\partial_{h}^{+j}f(mh) & := & h^{-1}(f(mh+e_jh)-f(mh)), \\
\partial_{h}^{-j}f(mh) & := & h^{-1}(f(mh)-f(mh-e_jh)),
\end{array}
\end{equation}
for discrete functions $f(mh)$ with $mh\in h\bbbz^{n}$. In the sequel, we consider functions defined on $\Omega_{h} \subset  h\bbbz^{8}$ and taking values in octonions $\bbbo$. As usual, all important properties such as, $l^{p}$-summability ($1\leq p<\infty$), are defined component-wisely.\par
Next step is to introduce discretisation of the Cauchy-Riemann operators in octonions. Several approaches to the discretisation of the Cauchy-Riemann (and Dirac) operators have been presented in recent years. In particular, the discrete Clifford analysis is generally based on the idea of splitting each basis element $e_{k}, k=0,1,\ldots,7$, into two basis elements $e_{k}^{+}$ and $e_{k}^{-}, k=0,1,\ldots,7$, i.e., $e_{k}=e_k^{+}+e_k^{-}$, corresponding to the forward and backward directions, respectively, see \cite{Brackx_1,FKS07} for the details. A typical choice for such a basis is one satisfying the relations:
\begin{eqnarray*}
\left\{
\begin{array}{ccc}
e_j^{-}e_k^{-}+e_k^{-}e_j^{-} &=&0, \\
e_j^{+}e_k^{+}+e_k^{+}e_j^{+}&=&0, \\
e_j^{+}e_k^{-} + e_k^{-}e_j^{+}&=&-\delta_{jk},
\end{array}
\right.
\end{eqnarray*}
where $\delta_{jk}$ is the Kronecker delta. This approach has several advantages and, in particular, it leads to a canonical factorisation of a star-Laplacian $\Delta_{h}$ by a pair of discrete Dirac operators. Unfortunately, this approach is not so well suited for working in the octonionic setting, because it is not so easy to respect the non-associativity of octonionic multiplication.\par
Another way of working with discrete Cauchy-Riemann and Dirac operators is to represent these operators by help of matrices containing finite difference approximations of partial derivatives, see for example \cite{Cerejeiras_1,FGHK06,Guerlebeck_1} and references therein. Similar to the first approach, using matrix-based discretisation for discretising the octonionic analysis will be difficult because of non-associativity, which is not respected by the classical matrix multiplication.\par
For proposing a discretisation of octonionic analysis respecting the non-associativity of octonionic multiplication, we will work with the approach presented in \cite{Faustino_1} and consisting in a direct discretisation of the continuous Dirac operators by forward and backward finite difference operators. Thus, by help of the finite difference operators~(\ref{Finite_differences}), we introduce {\itshape discrete forward Cauchy-Riemann operator} $D^{+}\colon l^{p}(\Omega_{h},\bbbo)\to l^{p}(\Omega_{h},\bbbo)$ and {\itshape discrete backward Cauchy-Riemann operators} $D^{-}\colon l^{p}(\Omega_{h},\bbbo)\to l^{p}(\Omega_{h},\bbbo)$ as follows
\begin{equation}
\label{Cauchy_Riemann_operators_discrete}
D^{+}_{h}:=\sum_{j=0}^{7} e_j\partial_{h}^{+j}, \quad D^{-}_{h}:=\sum_{j=0}^{7} e_j\partial_{h}^{-j}.
\end{equation}
A small disadvantage of this approach is related to the factorisation of the star-Laplacian, which is not just a composition of the Cauchy-Riemann operator and its conjugated operator, but requires a more complicated combination. It is easy to show by direct computations, that the star-Laplacian $\Delta_{h}$ can be represented as follows:
\begin{eqnarray*}
\Delta_{h} = \frac{1}{2}\left(D_{h}^{+}\overline{D_{h}^{-}}+D_{h}^{-}\overline{D_{h}^{+}}\right) \mbox{ with } \Delta_{h}:=\sum_{j=1}^{n}\partial_{h}^{+j}\partial_h^{-j},
\end{eqnarray*}
where $\overline{D_{h}^{-}}$ and $\overline{D_{h}^{+}}$ represent conjugated operators:
\begin{eqnarray*}
\overline{D_{h}^{-}}=\partial_{h}^{-0}-\sum_{j=1}^{7} e_j\partial_{h}^{-j}, \quad \overline{D_{h}^{+}}=\partial_{h}^{+0}-\sum_{j=1}^{7} e_j\partial_{h}^{+j}.
\end{eqnarray*}\par
In the rest of the paper, we will work with the discrete Cauchy-Riemann operators~(\ref{Cauchy_Riemann_operators_discrete}), because this discretisation clearly respects the non-associativity of octonionic multiplication.\par

\section{The discrete Stokes' formula in octonions}

In this section, we introduce the discrete Stokes' formula in octonionic setting. Additionally, we will underline the difference between octonionic constructions and the classical discrete Clifford analysis. Moreover, for keeping notations shorter, we will omit the lattice constant $h$ in the argument of discrete functions for the proof of discrete Stokes' formula, i.e. notations $f(m)$ or $f(m_{1},m_{2},m_{3})$ will be used instead of $f(mh)$ or $f(m_{1}h,m_{2}h,m_{3}h)$, respectively.\par
The following theorem presents the discrete octonionic Stokes' formula for the whole space:
\begin{theorem} 
The discrete Stokes' formula for the whole space with the lattice $h\bbbz^{8}$ is given by
\begin{equation}
\label{DiscreteStokesFormula}
\begin{array}{rl}
& \displaystyle \sum_{m\in \bbbz^{8}}  \left\{ \left[ g(mh)D_h^{+}\right] f(mh) - g(mh) \left[ D_h^{-}f(mh) \right]  \right\} h^8 = 0
\end{array}  
\end{equation}
for all discrete functions $f$ and $g$ such that the series converge.
\end{theorem}
\begin{proof} 
To underline clearly the effect of non-associativity of the octonionic multiplication, the proof will be presented with all explicit calculations. We start the proof by working with the first term on the left-hand side in~(\ref{DiscreteStokesFormula}):
\begin{eqnarray*}
\begin{array}{rl}
& \displaystyle \sum\limits_{m\in \bbbz^{8}} \left[g(m)D_{h}^{+}\right]f(m)h^{8} = \sum\limits_{m\in \bbbz^{8}} \sum\limits_{j=0}^{7} \left[\partial^{+j}g(m)e_{j}\right] f(m) h^{8} \\
\\
= & \displaystyle \sum\limits_{m\in \bbbz^{8}} \sum\limits_{j=0}^{7}\sum \limits_{i=0}^{7}\sum \limits_{k=0}^{7} \left[\partial^{+j}g_{i}(m)e_{i}e_{j}\right] f_{k}(m)e_{k} h^{8}.
\end{array}
\end{eqnarray*}
Next, using the relation $(e_{i}e_{j})e_{k}=-e_{i}(e_{j}e_{k})$ and the definition of $D_{h}^{+}$ leads to the following expression
\begin{eqnarray*}
\begin{array}{rl}
& \displaystyle \sum\limits_{m\in \bbbz^{8}} \sum\limits_{j=0}^{7}\sum \limits_{i=0}^{7}\sum \limits_{k=0}^{7} \left[-\partial^{+j}g_{i}f_{k}(m)e_{i}(e_{j}e_{k})\right]h^{8} \\
\\
= & \displaystyle \sum\limits_{m\in \bbbz^{8}} \sum\limits_{j=0}^{7}\sum \limits_{i=0}^{7}\sum \limits_{k=0}^{7} \left[-\left(g_{i}(m+e_{j})-g_{i}(m)\right)f_{k}(m)e_{i}(e_{j}e_{k})\right]h^{8} \\
\\
= & \displaystyle \sum\limits_{m\in \bbbz^{8}} \sum\limits_{j=0}^{7}\sum \limits_{i=0}^{7}\sum \limits_{k=0}^{7} \left[-g_{i}(m+e_{j})e_{i}f_{k}(m)+g_{i}(m)e_{i}f_{k}(m)\right](e_{j}e_{k})h^{8}.
\end{array}
\end{eqnarray*}
Performing change of variables in the last expression, we get
\begin{eqnarray*}
\begin{array}{rl}
& \displaystyle \sum\limits_{m\in \bbbz^{8}} \sum\limits_{j=0}^{7}\sum \limits_{i=0}^{7}\sum \limits_{k=0}^{7} \left[-g_{i}(m)e_{i}f_{k}(m-e_{j})+g_{i}(m)e_{i}f_{k}(m)\right](e_{j}e_{k})h^{8} \\
\\
= & \sum\limits_{m\in \bbbz^{8}} \sum\limits_{j=0}^{7}\sum \limits_{i=0}^{7}\sum \limits_{k=0}^{7} \left[g_{i}(m)e_{i}\left(f_{k}(m-e_{j})+f_{k}(m)\right)\right](e_{j}e_{k})h^{8} \\
\\
= & \sum\limits_{m\in \bbbz^{8}} \sum\limits_{j=0}^{7}\sum \limits_{i=0}^{7}\sum \limits_{k=0}^{7} g_{i}(m)e_{i}\partial^{-j} f_{k}(e_{j}e_{k})h^{8} \\
\\
= & \sum\limits_{m\in \bbbz^{8}} \sum\limits_{j=0}^{7}\sum \limits_{i=0}^{7}\sum \limits_{k=0}^{7} g_{i}(m)e_{i}\left(\partial^{-j}e_{j} f_{k}e_{k}\right)h^{8} = \sum\limits_{m\in \bbbz^{8}} g(m)\left[D_{h}^{-}f(m)\right]h^{8}.
\end{array}
\end{eqnarray*}
Thus, the statement of the theorem is proved.
\end{proof}\par
As we see from this theorem, the discrete Stokes' formula does not contain the associator in contrast to the continuous case~(\ref{Stokes_formula_continuous}). This is an interesting result, and a possible reason for vanishing of the associator could be the fact, that the discrete octonionic Stokes' formula contains two different differential operators: forward and backward Cauchy-Riemann operators, while in the continuous case both operators are the same. Additionally, it is worth to underline that the non-associativity affect the sign of the second summand in~(\ref{DiscreteStokesFormula}), which is not the case in the discrete Clifford analysis \cite{CKKS,Cerejeiras}.\par
Next, we consider the case of the upper half-lattice, defined as follows
\begin{eqnarray*}
h\bbbz_{+}^{8} := \left\{(h\underline{m},hm_{7})\colon \underline{m}\in\bbbz^{7},m_{7}\in\bbbz_{+}\right\}.
\end{eqnarray*}
The discrete octonionic Stokes' formula for the upper half-lattice is provided by the following theorem:
\begin{theorem} 
The discrete Stokes' formula for the upper half-lattice $h\bbbz_{+}^{8}$ is given by
\begin{equation}
\label{DiscreteStokesFormula_half_lattice}
\begin{array}{c}
\displaystyle \sum_{m\in \bbbz_{+}^{8}}  \left\{ \left[ g(mh)D_h^{+}\right] f(mh) - g(mh) \left[ D_h^{-}f(mh) \right]  \right\} h^8 \\
\\
\displaystyle = \sum\limits_{\underline{m}\in \bbbz^{7}} e_{7}\left(g(\underline{m},1)f_{k}(\underline{m},0)\right) h^{8}
\end{array}  
\end{equation}
for all discrete functions $f$ and $g$ such that the series converge.
\end{theorem}
\begin{proof} 
The proof of this theorem is similar to the proof of the discrete Stokes' formula for the whole space. Nonetheless, it is necessary to address the fact, that the discrete Cauchy-Riemann operators can be applied only for points with $m_{7}\geq 1$.  We start the proof by working with the first term on the left-hand side in~(\ref{DiscreteStokesFormula_half_lattice}):
\begin{eqnarray*}
\begin{array}{rl}
& \displaystyle \sum\limits_{m\in \bbbz_{+}^{8}} \left[g(m)D_{h}^{+}\right]f(m)h^{8} = \sum\limits_{m\in \bbbz_{+}^{8}} \sum\limits_{j=0}^{6} \left[\partial^{+j}g(m)e_{j}\right] f(m) h^{8} \\
\\
& \displaystyle + \sum\limits_{m\in \bbbz_{+}^{8}} \left[\partial^{+7}g(m)e_{7}\right] f(m) h^{8} = \displaystyle \sum\limits_{m\in \bbbz^{8}_{+}} \sum\limits_{j=0}^{6}\sum \limits_{i=0}^{7}\sum \limits_{k=0}^{7} \left[\partial^{+j}g_{i}(m)e_{i}e_{j}\right] f_{k}(m)e_{k} h^{8} \\
\\
& \displaystyle + \sum\limits_{\underline{m}\in \bbbz^{7}} \left\{\sum\limits_{m_{7}\geq 1}\sum \limits_{i=0}^{7}\sum \limits_{k=0}^{7} \left[\left(g_{i}(m+e_{7})f_{k}(m) - g_{i}(m)f_{k}(m)\right)e_{7}e_{i}\right] e_{k} h^{8} \right\}.
\end{array}
\end{eqnarray*}
Next, we will work with the second sum. By using the relation $(e_{i}e_{j})e_{k}=-e_{i}(e_{j}e_{k})$ and performing change of variables, we get the following expression
\begin{eqnarray*}
\begin{array}{rl}
& \displaystyle \sum\limits_{\underline{m}\in \bbbz^{7}} \left\{\sum\limits_{m_{7}\geq 1}\sum \limits_{i=0}^{7}\sum \limits_{k=0}^{7} \left[\left(-g_{i}(m+e_{7})f_{k}(m) + g_{i}(m)f_{k}(m)\right)e_{7}\right] e_{i}e_{k} h^{8} \right\} \\
\\
= & \displaystyle \sum\limits_{\underline{m}\in \bbbz^{7}} \left\{\sum\limits_{m_{7}\geq 1}\sum \limits_{i=0}^{7}\sum \limits_{k=0}^{7} \left[\left(g_{i}(m)f_{k}(m)-g_{i}(m)f_{k}(m-e_{7})\right)e_{7}\right] e_{i}e_{k} h^{8} \right\} \\
\\
= & \displaystyle \sum\limits_{\underline{m}\in \bbbz^{7}} \left\{\sum\limits_{m_{7}\geq 1}\sum \limits_{i=0}^{7}\sum \limits_{k=0}^{7} g_{i}(m)f_{k}(m)e_{7}\left(e_{i}e_{k}\right) h^{8} \right. \\
\\
& - \displaystyle \left.\sum\limits_{m_{7}\geq 2}\sum \limits_{i=0}^{7}\sum \limits_{k=0}^{7} g_{i}(m)f_{k}(m-e_{7})e_{7}\left(e_{i}e_{k}\right) h^{8} \right\} \\
\\
= & \displaystyle \sum\limits_{\underline{m}\in \bbbz^{7}} \left\{\sum\limits_{m_{7}\geq 1}\sum \limits_{i=0}^{7}\sum \limits_{k=0}^{7} g_{i}(m)f_{k}(m)e_{7}\left(e_{i}e_{k}\right) h^{8} \right. \\
\\
& - \displaystyle \sum\limits_{m_{7}\geq 1}\sum \limits_{i=0}^{7}\sum \limits_{k=0}^{7} g_{i}(m)f_{k}(m-e_{7})e_{7}\left(e_{i}e_{k}\right) h^{8} \\
\\
& + \displaystyle \left. \sum \limits_{i=0}^{7}\sum \limits_{k=0}^{7} g_{i}(\underline{m},1)f_{k}(\underline{m},0)e_{7}\left(e_{i}e_{k}\right) h^{8}\right\}.
\end{array}
\end{eqnarray*}
Combining this result with the first sum of the original expression, we finally get the following equality
\begin{eqnarray*}
\begin{array}{rcl}
\displaystyle \sum\limits_{m\in \bbbz_{+}^{8}} \left[g(m)D_{h}^{+}\right]f(m)h^{8} & = & \displaystyle \sum\limits_{m\in \bbbz_{+}^{8}} g(m)\left[D_{h}^{-}f(m)\right]h^{8}\\
\\
& & \displaystyle + \sum\limits_{\underline{m}\in \bbbz^{7}} e_{7}\left(g(\underline{m},1)f_{k}(\underline{m},0)\right) h^{8},
\end{array}
\end{eqnarray*}
which proofs the assertion of the theorem.
\end{proof}\par
Similarly, a discrete Stokes' formula can be established for the lower half-lattice, defined as follows
\begin{eqnarray*}
h\bbbz_{-}^{8} := \left\{(h\underline{m},hm_{7})\colon \underline{m}\in\bbbz^{7},m_{7}\in\bbbz_{-}\right\}.
\end{eqnarray*}
We have then the following theorem:
\begin{theorem} 
The discrete Stokes' formula for the lower half-lattice $h\bbbz_{-}^{8}$ is given by
\begin{equation}
\begin{array}{c}
\displaystyle \sum_{m\in \bbbz_{-}^{8}}  \left\{ \left[ g(mh)D_h^{+}\right] f(mh) - g(mh) \left[ D_h^{-}f(mh) \right]  \right\} h^8 \\
\\
\displaystyle = -\sum\limits_{\underline{m}\in \bbbz^{7}} e_{7}\left(g(\underline{m},0)f_{k}(\underline{m},-1)\right) h^{8}
\end{array}  
\end{equation}
for all discrete functions $f$ and $g$ such that the series converge.
\end{theorem}
\begin{proof} 
The proof of this theorem is analogue to the previous proof, and, therefore, we will present a shorter version of the proof. Hence, we have:
\begin{eqnarray*}
\begin{array}{rl}
& \displaystyle \sum\limits_{m\in \bbbz_{-}^{8}} \left[g(m)D_{h}^{+}\right]f(m)h^{8} = \sum\limits_{m\in \bbbz^{8}_{-}} \sum\limits_{j=0}^{6}\sum \limits_{i=0}^{7}\sum \limits_{k=0}^{7} \left[\partial^{+j}g_{i}(m)e_{i}e_{j}\right] f_{k}(m)e_{k} h^{8} \\
\\
& \displaystyle + \sum\limits_{\underline{m}\in \bbbz^{7}} \left\{\sum\limits_{m_{7}\leq -1}\sum \limits_{i=0}^{7}\sum \limits_{k=0}^{7} \left[\left(g_{i}(m+e_{7})f_{k}(m) - g_{i}(m)f_{k}(m)\right)e_{7}e_{i}\right] e_{k} h^{8} \right\}.
\end{array}
\end{eqnarray*}
Working with the second sum, we get
\begin{eqnarray*}
\begin{array}{rl}
& \displaystyle \sum\limits_{\underline{m}\in \bbbz^{7}} \left\{\sum\limits_{m_{7}\leq -1}\sum \limits_{i=0}^{7}\sum \limits_{k=0}^{7} \left[\left(g_{i}(m+e_{7})f_{k}(m) - g_{i}(m)f_{k}(m)\right)e_{7}e_{i}\right] e_{k} h^{8} \right\} \\
\\
= & \displaystyle \sum\limits_{\underline{m}\in \bbbz^{7}} \left\{\sum\limits_{m_{7}\leq -1}\sum \limits_{i=0}^{7}\sum \limits_{k=0}^{7} g_{i}(m)f_{k}(m)e_{7}\left(e_{i}e_{k}\right) h^{8} \right. \\
\\
& - \displaystyle \left.\sum\limits_{m_{7}\leq 0}\sum \limits_{i=0}^{7}\sum \limits_{k=0}^{7} g_{i}(m)f_{k}(m-e_{7})e_{7}\left(e_{i}e_{k}\right) h^{8} \right\} \\
\\
= & \displaystyle \sum\limits_{\underline{m}\in \bbbz^{7}} \left\{\sum\limits_{m_{7}\leq -1}\sum \limits_{i=0}^{7}\sum \limits_{k=0}^{7} g_{i}(m)f_{k}(m)e_{7}\left(e_{i}e_{k}\right) h^{8} \right. \\
\\
& - \displaystyle \sum\limits_{m_{7}\leq -1}\sum \limits_{i=0}^{7}\sum \limits_{k=0}^{7} g_{i}(m)f_{k}(m-e_{7})e_{7}\left(e_{i}e_{k}\right) h^{8} \\
\\
& - \displaystyle \left. \sum \limits_{i=0}^{7}\sum \limits_{k=0}^{7} g_{i}(\underline{m},0)f_{k}(\underline{m},-1)e_{7}\left(e_{i}e_{k}\right) h^{8}\right\}.
\end{array}
\end{eqnarray*}
Combining this result with the first sum of the original expression, we obtain the assertion of the theorem.
\end{proof}\par

\section{Summary}

While a lot of results in the continuous octonionic analysis have been presented in recent years, construction of a discrete counterpart of the continuous theory is still missing. Therefore, in this short paper, we discussed first ideas towards developing a discrete octonionic analysis. In particular, we discuss several approaches to a discretisation of octonionic analysis, and underlined, that because of non-associativity of octonions not all approaches common in the discrete Clifford analysis are applicable in the octonionic setting. After that, we presented several discrete octonionic Stokes' formulae: for the whole spaces, upper-half lattice, and lower-half lattice.  The results presented in this paper will be further extended in future work.


\end{document}